\documentclass[10pt,twoside,reqno]{amsart}
\usepackage{amssymb}
\usepackage{multirow}
\usepackage{graphicx}
\usepackage[dvips]{epsfig}
\usepackage{latexsym}
\usepackage{amsmath}
\usepackage{amsthm}
\usepackage{amscd}
\usepackage{amssymb}
\input xy
\xyoption{all}

\calclayout

\numberwithin{equation}{section}
\makeatletter

\renewcommand{\@secnumfont}{\bfseries}

\renewcommand{\section}{\@startsection{section}{1}%
  {0mm}{.7\linespacing\@plus\linespacing}{.5\linespacing}
  {\normalfont\bfseries\centering}}

\newcommand{\bibsection}{\@startsection{section}{1}%
  {0mm}{.7\linespacing\@plus\linespacing}{.5\linespacing}
  {\normalfont\scshape\centering}}

\renewcommand{\@biblabel}[1]{#1.}

\newtheorem{thm}{\bf Theorem}[section]

\theoremstyle{theorem}

\theoremstyle{definition}

\numberwithin{equation}{section}

\begin{document}

\vspace{1.3cm}

\title {A new type degenerate Daehee numbers and polynomials}

\author{Taekyun Kim}
\address{Department of Mathematics, kwangwoon University, Seoul 139-701, Republic of Korea}
\email{tkkim@kw.ac.kr}

\author{Dae San Kim}
\address{Department of Mathematics, Sogang University, Seoul 121-742, Republic of Korea}
\email{dskim@sogang.ac.kr}

\author{Han Young kim}
\address{Department of Mathematics, kwangwoon University, Seoul 139-701, Republic of Korea}
\email{gksdud213@kw.ac.kr}

\author{Jongkyum Kwon}
\address{Department of Mathematics Education and ERI, Gyeongsang National University, Jinju, Gyeongsangnamdo 52828, Korea}
\email{mathkjk26@gnu.ac.kr}

\thanks{}

\begin{abstract}
In this paper, we study a degenerate version of the Daehee polynomials and numbers, namely the degenerate Daehee polynomials and numbers, which were recently introduced by Jang et. al. We derive their explicit expressions and some identities involving them. Further, we introduce the multiple degenerate Daehee numbers and higher-order degenerate Daehee polynomials and numbers which can be represented in terms of integrals on the unitcube. Again, we deduce their explicit expressions and some identities related to them.
\end{abstract}

\maketitle

\bigskip
\medskip
\section{\bf Introduction}
\medskip
The degenerate versions of Bernoulli and Euler polynomials, namely the degenerate Bernoulli and Euler polynomials, were studied by Carlitz in \cite{1}. In recent years, studying various degenerate versions of some special polynomials and numbers drew attention of some mathematicians and many arithmetic and combinatorial results were obtained [4,5,7-9,11-14,17]. They can be explored by using various tools like combinatorial methods, generating functions, differential equations, umbral calculus techniques, $p$-adic analysis, and probability theory. \\
\indent The aim of this paper is to study a degenerate version of the Daehee polynomials and numbers, namely the degenerate Daehee polynomials and numbers, in the spirit of \cite{1}. They were recently introduced by Jang et. al \cite{4}. We derive their explicit expressions and some identities involving them. Further, we introduce the multiple degenerate Daehee numbers and higher-order degenerate Daehee polynomials and numbers. Again, we deduce their explicit expressions and some identities related to them. \\
\indent This paper is organized as follows. In Section 1, we state what we need in the rest of the paper. These include the Stirling numbers of the first and second kinds, the higher-order Bernoulli polynomials, the higher-order Daehee polynomials, the higher-order degenerate Bernoulli polynomials, the degenerate exponential functions, and the degenerate Stirling numbers of the first and second kinds. In Section 2, we define the degenerate Daehee polynomials and numbers whose generating functions can be expressed in terms of integrals on the unit interval.  We find their explicit expressions and some identities involving them. We also introduce the multiple degenerate Daehee numbers, the generating function of which can be expressed in terms of a multiple integral on the unitcube or of the modified polyexponential function \cite{6}. We deduce an explicit expression of them and some identities involving them. In Section 3, we introduce the higher-order degenerate Daehee polynomials and numbers whose generating function can be represented as a multiple integral on the unitcube. We derive their explicit expressions and some identities relating to them. Finally, we conclude this paper in Section 4.

\medskip

For $n \ge 0$, the Stirling numbers of the first kind are defined by
\begin{equation}\label{01}
\begin{split}
(x)_{n} = \sum_{l=0}^{n}S_{1}(n,l)x^{l} ,  \,\,(\mathrm{see}\,\, [1,3]),
\end{split}
\end{equation}
where $(x)_{0} = 1$, $(x)_{n} = x(x-1)\cdots (x-n+1),\, (n \geq 1).$\\
As an inversion formula of \eqref{01}, the Stirling numbers of the second kind are defined as
\begin{equation}\label{02}
\begin{split}
x^{n} = \sum_{l=0}^{n}S_{2}(n,l)(x)_{l}, (n \ge 0), \,\,(\mathrm{see}\,\, [3]).
\end{split}
\end{equation}
For $\alpha \in \mathbb{N}$, the Bernoulli polynomials of order $\alpha$ are defined as
\begin{equation}\label{03}
\begin{split}
\left(\frac{t}{e^{t}-1}\right)^{\alpha}e^{xt} = \sum_{n=0}^{\infty}B_{n}^{(\alpha)}(x)\frac{t^{n}}{n!}, \,\,(\mathrm{see}\,\, [3]).
\end{split}
\end{equation}
$B_{n}(x) = B_{n}^{(1)}(x)$ are called the Bernoulli polynomials and  $B^{(\alpha)}_{n} = B^{(\alpha)}_{n}(0)$ the Bernoulli numbers of order $\alpha$.\\
The Daehee polynomials are defined by
\begin{equation}\label{04}
\begin{split}
\left( \frac{\log(1+t)}{t}\right)(1+t)^{x} = \sum_{n=0}^{\infty}D_{n}(x)\frac{t^{n}}{n!},\,\, (\mathrm{see}\,\, [2,3,10,15-34]).
\end{split}
\end{equation}
For $x=0$, $D_{n} = D_{n}(0)$ are called the Daehee numbers.\\
Recently, EI-Desouky and Mustafa studied some new results on the higher--order Daehee and Bernoulli numbers and polynomials (see [3]).\\
That is, they derived a new matrix representation for the higher--order Daehee number and polynomials, the higher--order $\lambda$--Daehee numbers and polynomials, and the higher-order twisted $\lambda$--Daehee numbers and polynomials (see [3,4]).\\
The Daehee polynomials of order $k$ are defined by
\begin{equation}\label{05}
\begin{split}
\left(\frac{\log(1+t)}{t}\right)^{k}(1+t)^{x} = \sum_{n=0}^{\infty}D_{n}^{(k)}(x)\frac{t^{n}}{n!}, \,\, (\mathrm{see}\,\, [3,20]).
\end{split}
\end{equation}
In \cite{3}, EI--Desouky and Mustafu obtained the following interesting formula.
\begin{equation}\label{06}
\begin{split}
D_{m}^{(k)}(z) = m!\sum_{n=0}^{m}{z \choose m-n}b_{m}^{(-k)},
\end{split}
\end{equation}
where $b_{n}^{(x)}$ are the N\"orlund numbers of the second kind given by (see [25])
\begin{equation}\label{06-1}
\bigg(\frac{t}{\log (1+t)}\bigg)^x=\sum_{n=0}^{\infty}b_n^{(x)}t^n.
\end{equation}

Recently, Daehee numbers and polynomials have been studied by many researchers in various areas (see  [2,3,10,15-34]).\\
In \cite{1}, Carlitz considered the degenerate Bernoulli polynomials given by
\begin{equation}\label{07}
\begin{split}
\frac{t}{(1+\lambda t)^{\frac{1}{\lambda}}-1}(1+\lambda t)^{\frac{x}{\lambda}} = \sum_{n=0}^{\infty}\beta_{n,\lambda}^{(x)}\frac{t^{n}}{n!}, \,(\lambda \in \mathbb{R}).
\end{split}
\end{equation}
When $x=0$, $\beta_{n,\lambda} = \beta_{n,\lambda}(0)$ are called the degenerate Bernoulli numbers.
For $r \in \mathbb{N}$, he also defined the higher--order degenerate Bernoulli polynomials as
\begin{equation}\label{08}
\begin{split}
\left(\frac{t}{(1+\lambda t)^{\frac{1}{\lambda}}-1}\right)^{r}(1+\lambda t)^{\frac{x}{\lambda}} = \sum_{n=0}^{\infty}\beta_{n,\lambda}^{(r)}(x)\frac{t^{n}}{n!}, \,\,(\mathrm{see}\,\, [7]).
\end{split}
\end{equation}
When $x=0$, $\beta_{n,\lambda}^{(r)} = \beta_{n,\lambda}^{(r)}(0)$ are called the degenerate Bernoulli numbers of order $r$.\\
The degenerate exponential functions are given by
\begin{equation}\label{09}
\begin{split}
e^{x}_{\lambda}(t) = (1+\lambda t)^{\frac{x}{\lambda}}, e_{\lambda}(t)  = e^{1}_{\lambda}(t) = (1+\lambda t)^{\frac{1}{\lambda}}, \,\,(\mathrm{see}\,\, [3,11]).
\end{split}
\end{equation}
We note that
\begin{equation}\label{10}
\begin{split}
e^{x}_{\lambda}(t) = \sum_{n=0}^{\infty}\frac{(x)_{n,\lambda}}{n!}t^{n}, \,\,(\mathrm{see}\,\, [11]),
\end{split}
\end{equation}
where $(x)_{0,\lambda} = 1$, $(x)_{n,\lambda} = x(x-\lambda)\cdots (x-(n-1)\lambda)$, $(n \ge 1)$.\\
Note that $\lim_{\lambda \to 0} e_{\lambda}^{x}(t) = e^{xt}$, $\lim_{\lambda \to  0}\beta^{(r)}_{n,\lambda}(x) = B_{n}^{(r)}(x)$.\\
Recently, Kim considered the degenerate Stirling numbers of the second kind given by
\begin{equation}\label{11}
\begin{split}
(x)_{n,\lambda} = \sum_{l=0}^{n}S_{2,\lambda}(n,l)(x)_{l}, (n \ge 0), \,\,(\mathrm{see}\,\, [11]).
\end{split}
\end{equation}
Note that $\lim_{\lambda \to 0}S_{2,\lambda}(n,l) = S_{2}(n,l)$.\\
From \eqref{11}, we note that
\begin{equation}\label{12}
\begin{split}
\frac{1}{k!}(e_{\lambda}(t)-1)^{k} = \sum_{n=k}^{\infty}S_{2,\lambda}(n,k)\frac{t^{n}}{n!}, \,(k \ge 0), \,\,(\mathrm{see}\,\, [11]).
\end{split}
\end{equation}
As an inversion formula of \eqref{11}, the Stirling numbers of the first kind are defined by
\begin{equation}\label{13}
\begin{split}
(x)_{n} = \sum_{l=0}^{n}S_{1,\lambda}(n,l)(x)_{l,\lambda},\,
(n \ge 0), \,\,(\mathrm{see}\,\, [11]).
\end{split}
\end{equation}
We see that $\log_{\lambda}(t)=\frac{1}{\lambda}(t^{\lambda}-1)$ is the compositional inverse of $e_{\lambda}(t)$ satisfying $\log_{\lambda}(e_{\lambda}(t))  = e_{\lambda}(\log_{\lambda}(t)) = t$.\\
By \eqref{13}, we get
\begin{equation}\label{14}
\begin{split}
\frac{1}{k!}\left(\log_{\lambda}(1+t)\right)^{k} = \sum_{n=k}^{\infty}S_{1,\lambda}(n,k)\frac{t^{n}}{n!}, \,\, (\mathrm{see}\,\, [11]).
\end{split}
\end{equation}
Note that $\lim_{\lambda \to 0}\log_{\lambda}(1+t) = \log(1+t)$.

\medskip

\section{\bf Degenerate Daehee numbers and polynomials}
\medskip

We define the degenerate Daehee polynomials by
\begin{equation}\label{15}
\begin{split}
\frac{\log_{\lambda}(1+t)}{t}(1+t)^{x} = \sum_{n=0}^{\infty}D_{n,\lambda}(x)\frac{t^{n}}{n!}, \,(\lambda \in \mathbb{R}).
\end{split}
\end{equation}
When $x=0$, $D_{n,\lambda} = D_{n,\lambda}(0)$ are called the degenerate Daehee numbers.\\
From \eqref{04} and \eqref{15}, we note that $\lim_{\lambda \to 0}D_{n,\lambda}(x) = D_{n}(x)$, $(n \ge 0)$.\\
We observe that
\begin{equation*}
\begin{split}
\frac{\log(1+t)}{t}\int_{0}^{1}(1+t)^{\lambda y+x}dy  = \frac{\log_{\lambda}(1+t)}{t}(1+t)^{x} = \sum_{n=0}^{\infty}D_{n,\lambda}(x)\frac{t^{n}}{n!}.
\end{split}
\end{equation*}
When $x=0$, we have
\begin{equation}\label{16}
\begin{split}
\frac{\log(1+t)}{t}\int_{0}^{1}(1+t)^{\lambda y} dy = \sum_{n=0}^{\infty}D_{n,\lambda}\frac{t^{n}}{n!}.
\end{split}
\end{equation}
On the other hand,
\begin{equation}\label{17}
\begin{split}
& \frac{\log(1+t)}{t}\int_{0}^{1}(1+t)^{\lambda y}dy  = \frac{\log(1+t)}{t}\sum_{m=0}^{\infty}\frac{\lambda^{m}(\log(1+t))^{m}}{(m+1)!}\\
& \quad  = \frac{1}{t}\sum_{m=0}^{\infty}\frac{(\log(1+t))^{m+1}}{(m+1)!}\lambda^{m} = \frac{1}{t}\sum_{m=1}^{\infty}\lambda^{m-1}\frac{1}{m!}(\log(1+t))^{m} \\
&\quad = \frac{1}{t}\sum_{m=1}^{\infty}\lambda^{m-1}\sum_{n=m}^{\infty}S_{1}(n,m)\frac{t^{n}}{n!} = \frac{1}{t}\sum_{n=1}^{\infty}\left(\sum_{m=1}^{n}\lambda^{m-1}S_{1}(n,m)\right)\frac{t^{n}}{n!}\\
& \quad = \sum_{n=0}^{\infty}\left(\frac{1}{n+1}\sum_{m=1}^{n+1}\lambda^{m-1}S_{1}(n+1,m)\right)\frac{t^{n}}{n!}.
\end{split}
\end{equation}

Therefore, by \eqref{16} and \eqref{17}, we obtain the following theorem.

\begin{thm}\label{Theorem 1}
For $n\ge 0$, we have
\begin{equation*}
\begin{split}
D_{n,\lambda} = \frac{1}{n+1}\sum_{m=1}^{n+1}\lambda^{m-1}S_{1}(n+1,m).
\end{split}
\end{equation*}
\end{thm}
By replacing $t$ by $e_{\lambda}(t)-1$ in \eqref{15}, we get
\begin{equation}\label{18}
\begin{split}
\frac{t}{e_{\lambda}(t)-1}e_{\lambda}^{x}(t) & = \sum_{m=0}^{\infty}D_{m,\lambda}(x)\frac{1}{m!}(e_{\lambda}(t)-1)^{m}\\
& = \sum_{m=0}^{\infty}D_{m,\lambda}(x)\sum_{n=m}^{\infty}S_{2,\lambda}(n,m)\frac{t^{n}}{n!}\\
& = \sum_{n=0}^{\infty}\left(\sum_{m=0}^{n}D_{m,\lambda}(x)S_{2,\lambda}(n,m)\right)\frac{t^{n}}{n!}.
\end{split}
\end{equation}
On the other hand,
\begin{equation}\label{19}
\begin{split}
\frac{t}{e_{\lambda}(t)-1}e_{\lambda}^{x}(t) = \sum_{n=0}^{\infty}\beta_{n,\lambda}(x)\frac{t^{n}}{n}.
\end{split}
\end{equation}
Therefore, by \eqref{18} and \eqref{19}, we obtain the following theorem.
\begin{thm}\label{Theorem 2}
For $n \ge 0$, we have
\begin{equation*}
\begin{split}
\beta_{n,\lambda}(x) = \sum_{m=0}^{n} D_{m,\lambda}(x)S_{2,\lambda}(n,m).
\end{split}
\end{equation*}
\end{thm}
Note that
\begin{equation*}
\begin{split}
B_{n}(x) = \lim_{\lambda \to 0}\beta_{n,\lambda}(x) = \sum_{m=0}^{n}D_{m}(x)S_{2}(n,m), (n \ge 0).
\end{split}
\end{equation*}
To find the inversion formula of Theorem 2, we replace $t$ by $\log_{\lambda}(1+t)$ in \eqref{07} and get
\begin{equation}\label{20}
\begin{split}
\frac{\log_{\lambda}(1+t)}{t}(1+t)^{x}& = \sum_{m=0}^{\infty}\beta_{m,\lambda}(x)\frac{1}{m!}\left(\log_{\lambda}(1+t)\right)^{m}\\
& = \sum_{m=0}^{\infty}\beta_{m,\lambda}(x)\sum_{n=m}^{\infty}S_{1,\lambda}(n,m)\frac{t^{n}}{n!}\\
& = \sum_{n=0}^{\infty}\left(\sum_{m=0}^{n}\beta_{m,\lambda}(x)S_{1,\lambda}(n,m)\right)\frac{t^{n}}{n!}.
\end{split}
\end{equation}
Therefore, by \eqref{15} and \eqref{20}, we obtain the following theorem.
\begin{thm}\label{Theorem 3}
For $n \ge 0$, we have
\begin{equation*}
\begin{split}
D_{n,\lambda}(x) = \sum_{m=0}^{n}\beta_{m,\lambda}(x)S_{1,\lambda}(n,m).
\end{split}
\end{equation*}
\end{thm}
Note that
\begin{equation*}
\begin{split}
D_{n}(x) = \lim_{\lambda \to 0}D_{n,\lambda}(x) = \sum_{m=0}^{n}B_{m}(x)S_{1}(n,m),(n \ge 0).
\end{split}
\end{equation*}
From \eqref{09}, we can derive the following equation.
\begin{equation}\label{21}
\begin{split}
\sum_{n=0}^{\infty}D_{n,\lambda}(x)\frac{t^{n}}{n!} & = \frac{\log_{\lambda}(1+t)}{t}(1+t)^{x} = \frac{\log_{\lambda}(1+t)}{t}e_{\lambda}^{x}(\log_{\lambda}(1+t))\\
& = \frac{\log_{\lambda}(1+t)}{t} \sum_{m=0}^{\infty}(x)_{m,\lambda}\frac{(\log_{\lambda}(1+t))^{m}}{m!}\\
& = \frac{1}{t}\sum_{m=0}^{\infty}(m+1)(x)_{m,\lambda}\frac{1}{(m+1)!}\left(\log_{\lambda}(1+t)\right)^{m+1}\\
& = \frac{1}{t}\sum_{m=0}^{\infty}(m+1)(x)_{m,\lambda}\sum_{n=m+1}^{\infty}S_{1,\lambda}(n,m+1)\frac{t^{n}}{n!}\\
\end{split}
\end{equation}
\begin{equation}
\begin{split}
& = \sum_{m=0}^{\infty}(m+1)(x)_{m,\lambda}\sum_{n=m}^{\infty}\frac{S_{1,\lambda}(n+1,m+1)}{n+1}\frac{t^{n}}{n!}\\
& = \sum_{n=0}^{\infty}\left\{\frac{1}{n+1}\sum_{m=0}^{n}(m+1)(x)_{m,\lambda}S_{1,\lambda}(n+1,m+1)\right\}\frac{t^{n}}{n!}.
\end{split}
\end{equation}
Therefore, by \eqref{21}, we obtain the following theorem.
\begin{thm}\label{Theorem 4}
For $n \ge 0$, we have
\begin{equation*}
\begin{split}
D_{n,\lambda}(x) = \frac{1}{n+1}\sum_{m=0}^{n}(m+1)(x)_{m,\lambda}S_{1,\lambda}(n+1,m+1).
\end{split}
\end{equation*}
\end{thm}

For $s \in \mathbb{C}$, the polyexponential function is defined by Hardy as
\begin{equation}\label{22}
\begin{split}
e(x,a|s) = \sum_{n=0}^{\infty}\frac{x^{n}}{(n+a)^{s}n!}, (Re(a) > 0), \quad (see [7]).
\end{split}
\end{equation}
In [6], the modified polyexponential function is introduced as
\begin{equation}\label{23}
\begin{split}
Ei_{k}(x) = \sum_{n=1}^{\infty}\frac{x^{n}}{(n-1)!n^{k}}, (k \in \mathbb{Z}).
\end{split}
\end{equation}
Note that $x\,  e(x,1|k) = Ei_{k}(x)$.\\
We observe that
\begin{equation}\label{24}
\begin{split}
\frac{\partial}{\partial x_{1}}(1+t)^{\lambda x_{1}x_{2}\cdots x_{k}} = x_{2}\cdots x_{k}\lambda \log(1+t)(1+t)^{\lambda x_{1}x_{2}\cdots x_{k}}.
\end{split}
\end{equation}
Thus, by \eqref{24}, we get
\begin{equation}\label{25}
\begin{split}
&\frac{\log(1+t)}{t}\int_{0}^{1}(1+t)^{\lambda x_{1}x_{2}\cdots x_{k}}dx_{1} = \frac{1}{x_{2}x_{3}\cdots x_{k}}\frac{\log_{\lambda}(1+t)^{x_{2}\cdots x_{k}}}{t}\\
& \quad = \frac{1}{t}\frac{1}{x_{2}x_{3}\cdots x_{k}}\sum_{m=1}^{\infty}\lambda^{m-1}\frac{(\log(1+t))^{m}}{m!}x_{2}^{m}x_{3}^{m}\cdots x_{k}^{m}\\
& \quad = \frac{1}{t}\sum_{m=1}^{\infty}\frac{\lambda^{m-1}(\log(1+t))^{m}}{(m-1)!m} x_{2}^{m-1}x_{3}^{m-1}\cdots x_{k}^{m-1}.
\end{split}
\end{equation}
From \eqref{25}, we can derive the following equation:
\begin{equation}\label{26}
\begin{split}
& \frac{\log(1+t)}{t}\int_{0}^{1}\cdots\int_{0}^{1}(1+t)^{\lambda x_{1}x_{2}\cdots x_{k}}dx_{1}dx_{2}\cdots dx_{k}\\
& \quad = \frac{1}{t}\sum_{m=1}^{\infty}\frac{\lambda^{m-1}(\log(1+t))^{m}}{(m-1)!m^k} = \frac{1}{\lambda t}Ei_{k}(\lambda \log(1+t)).
\end{split}
\end{equation}
Now, we define the multiple degenerate Daehee numbers as the multiple integral on the unitcube given by
\begin{equation}\label{27}
\begin{split}
& \frac{\log(1+t)}{t}\int_{0}^{1}\cdots\int_{0}^{1}(1+t)^{\lambda x_{1}x_{2}\cdots x_{k}}dx_{1}dx_{2}\cdots dx_{k}\\
 & \quad = \sum_{n=0}^{\infty}\widehat{D}_{n,\lambda}^{(k)}\frac{t^{n}}{n!}.
\end{split}
\end{equation}
Then, by \eqref{26} and \eqref{27}, we get
\begin{equation}\label{28}
\begin{split}
\frac{1}{\lambda t}Ei_{k}(\lambda \log(1+t)) =  \sum_{n=0}^{\infty}\widehat{D}_{n,\lambda}^{(k)}\frac{t^{n}}{n!}.
\end{split}
\end{equation}
Note that $\widehat{D}_{n,\lambda}^{(1)} = D_{n,\lambda}$, $(
n \ge 0)$.\\
We observe that
\begin{equation}\label{29}
\begin{split}
& \frac{1}{\lambda t}Ei_{k}(\lambda \log(1+t)) = \frac{1}{\lambda t}\sum_{m=1}^{\infty}\frac{\lambda^{m}(\log(1+t))^{m}}{(m-1)!m^{k}}\\
&\quad = \frac{1}{\lambda t} \sum_{m=1}^{\infty}\frac{\lambda^{m}}{m^{k-1}}\frac{1}{m!}(\log(1+t))^{m} = \frac{1}{\lambda t}\sum_{m=1}^{\infty}\frac{\lambda^{m}}{m^{k-1}}\sum_{n=m}^{\infty}S_{1}(n,m)\frac{t^{n}}{n!}\\
& \quad =  \frac{1}{t}\sum_{n=1}^{\infty}\sum_{m=1}^{n}\frac{\lambda^{m-1}}{m^{k-1}} S_{1}(n,m)\frac{t^{n}}{n!}\\
& \quad = \sum_{n=0}^{\infty}\left(\frac{1}{n+1}\sum_{m=1}^{n+1}\frac{\lambda^{m-1}}{m^{k-1}}S_{1}(n+1,m)\right)\frac{t^{n}}{n!}.
\end{split}
\end{equation}
Therefore, by \eqref{28} and \eqref{29}, we obtain the following theorem.
\begin{thm}\label{Theorem 5}
For $n \ge 0$, we have
\begin{equation*}
\begin{split}
\widehat{D}_{n,\lambda}^{(k)} = \frac{1}{n+1}\sum_{m=1}^{n+1}\frac{\lambda^{m-1}}{m^{k-1}}S_{1}(n+1,m).
\end{split}
\end{equation*}
\end{thm}
By replacing $t$ by $e^{t}-1$ in \eqref{28}, we get
\begin{equation}\label{30}
\begin{split}
& \sum_{m=0}^{\infty}\widehat{D}_{m,\lambda}^{(k)}\frac{1}{m!}(e^{t}-1)^{m} = \frac{1}{\lambda(e^{t}-1)}Ei_{k}(\lambda t)\\
& \quad  = \frac{t}{e^{t}-1}\frac{1}{\lambda t}Ei_{k}(\lambda t) = \sum_{l=0}^{\infty}B_{l}\frac{t^{l}}{l!}\sum_{m=0}^{\infty}\frac{\lambda^{m}}{(m+1)^{k}}\frac{t^{m}}{m!}. \\
& \quad = \sum_{n=0}^{\infty}\left(\sum_{l=0}^{n}{n \choose l}\frac{\lambda^{n-l}B_{l}}{(n-l+1)^{k}}\right)\frac{t^{n}}{n!}.
\end{split}
\end{equation}
On the other hand,
\begin{equation}\label{31}
\begin{split}
&\sum_{m=0}^{\infty}\widehat{D}_{m,\lambda}^{(k)}\frac{1}{m!}(e^{t}-1)^{m}  = \sum_{m=0}^{\infty}\widehat{D}_{m,\lambda}^{(k)}\sum_{n=m}^{\infty}S_{2}(n,m)\frac{t^{n}}{n!}\\
& \quad = \sum_{n=0}^{\infty}\left(\sum_{m=0}^{n}\widehat{D}_{m,\lambda}^{(k)}S_{2}(n,m)\right)\frac{t^{n}}{n!}.
\end{split}
\end{equation}
Therefore, by \eqref{30} and \eqref{31}, we obtain the following theorem.
\begin{thm}\label{Theorem 6}
For $n \ge 0$, we have
\begin{equation*}
\begin{split}
\sum_{m=0}^{n}\widehat{D}_{m,\lambda}^{(k)}S_{2}(n,m) = \sum_{l=0}^{n}{n \choose l}\frac{\lambda^{n-l}B_{l}}{(n-l+1)^{k}}.\end{split}
\end{equation*}
\end{thm}
From \eqref{30}, we note that
\begin{equation}\label{32}
\begin{split}
\frac{1}{\lambda t}Ei_{k}(\lambda t) = & \frac{1}{t}(e^{t}-1)\sum_{m=0}^{\infty}\widehat{D}_{m,\lambda}^{(k)}\frac{1}{m!}(e^{t}-1)^{m}\\
& = \frac{1}{t}\sum_{m=1}^{\infty}m\widehat{D}_{m-1,\lambda}^{(k)}\frac{1}{m!}(e^{t}-1)^{m}\\
& = \frac{1}{t}\sum_{m=1}^{\infty}m \widehat{D}_{m-1,\lambda}^{(k)}\sum_{n=m}^{\infty}S_{2}(n,m)\frac{t^{n}}{n!}\\
&  = \frac{1}{t}\sum_{n=1}^{\infty}\sum_{m=1}^{n}m\widehat{D}_{m-1,\lambda}^{(k)}S_{2}(n,m)\frac{t^{n}}{n!}\\
& = \sum_{n=0}^{\infty}\left(\frac{1}{n+1}\sum_{m=1}^{n+1}m\widehat{D}_{m-1,\lambda}^{(k)}S_{2}(n+1,m)\right)\frac{t^{n}}{n!}.
\end{split}
\end{equation}
On the other hand,
\begin{equation}\label{33}
\begin{split}
\frac{1}{\lambda t}Ei_{k}(\lambda t) = \frac{1}{\lambda t}\sum_{n=1}^{\infty}\frac{\lambda^{n}t^{n}}{(n-1)!n^{k}} = \sum_{n=0}^{\infty}\frac{\lambda^{n}}{(n+1)^{k}}\frac{t^n}{n!}.
\end{split}
\end{equation}
Therefore, by \eqref{32} and \eqref{33}, we obtain the following theorem.
\begin{thm}\label{Theorem 7}
For $n \ge 0$, we have
\begin{equation*}
\begin{split}
\frac{\lambda^{n}}{(n+1)^{k}}  &= \frac{1}{n+1}\sum_{m=1}^{n}m\widehat{D}_{m-1,\lambda}^{(k)}S_{2}(n+1,m)\\
& = \frac{1}{n+1}\sum_{m=0}^{n-1}(m+1)\widehat{D}_{m,\lambda}^{(k)}S_{2}(n+1,m+1).
\end{split}
\end{equation*}
\end{thm}

\medskip

\section{\bf Higher-order degenerate Daehee numbers and polynomials}
\medskip

As an additive version of \eqref{27}, we consider the degenerate Daehee polynomials of order $r$ given by the following multiple integral on the unit cube\\
\begin{equation}\label{34}
\begin{split}
\sum_{n=0}^{\infty}D_{n,\lambda}^{(r)}(x)\frac{t^{n}}{n!} & = \left(\frac{\log(1+t)}{t}\right)^{r}\int_{0}^{1}\cdots \int_{0}^{1}(1+t)^{\lambda(x_{1}+\cdots +x_{r})+x}dx_{1}\cdots dx_{r}\\
& = \left(\frac{\log_{\lambda}(1+t)}{t}\right)^{r}(1+t)^{x}, (r\in \mathbb{N}).
\end{split}
\end{equation}
When $x=0$, $D_{n,\lambda}^{(r)} = D_{n,\lambda}^{(r)}(0)$, $(n \ge 0)$, are called the degenerate Daehee numbers of order $r$.\\
From \eqref{34}, we note that
\begin{equation}\label{35}
\begin{split}
\sum_{n=0}^{\infty}D_{n,\lambda}^{(r)}\frac{t^{n}}{n!} & =  \left(\frac{\log_{\lambda}(1+t)}{t}\right)^{r} = \frac{r!}{t^{r}}\frac{1}{r!}(\log_{\lambda}(1+t))^{r}\\
& = \frac{r!}{t^{r}}\sum_{n=r}^{\infty}S_{1,\lambda}(n,r)\frac{t^{n}}{n!}\\
& = \sum_{n=0}^{\infty}S_{1,\lambda}(n+r,r)\frac{r!n!}{(n+r)!}\frac{t^{n}}{n!}\\
& = \sum_{n=0}^{\infty}\frac{S_{1,\lambda}(n+r,r)}{{n+r \choose n}}\frac{t^{n}}{n!}.
\end{split}
\end{equation}
Therefore, by comparing the coefficients on both sides of \eqref{35}, we obtain the following theorem.
\begin{thm}\label{Theorem 8}
For $n \ge 0$, we have\begin{equation*}
\begin{split}
D_{n,\lambda}^{(r)} = \frac{1}{{n+r \choose n}}S_{1,\lambda}(n+r,r), (r \in \mathbb{N}).
\end{split}
\end{equation*}
\end{thm}
By replacing $t$ by $e_{\lambda}(t)-1$ in \eqref{34}, we get
\begin{equation}\label{36}
\begin{split}
\sum_{k=0}^{\infty}D_{k,\lambda}^{(r)}(x)\frac{1}{k!}(e_{\lambda}(t)-1)^{k} & = \left(\frac{t}{e_{\lambda}(t)-1}\right)^{r}e_{\lambda}^{x}(t)\\
& = \sum_{n=0}^{\infty}\beta_{n,\lambda}^{(r)}(x)\frac{t^{n}}{n!}.
\end{split}
\end{equation}
On the other hand,
\begin{equation}\label{37}
\begin{split}
&\sum_{k=0}^{\infty}D_{k,\lambda}^{(r)}(x)\frac{1}{k!}(e_{\lambda}(t)-1)^{k}\\
&\quad = \sum_{k=0}^{\infty}D_{k,\lambda}^{(r)}(x)\sum_{n=k}^{\infty}S_{2,\lambda}(n,k)\frac{t^{n}}{n!}\\
& \quad = \sum_{n=0}^{\infty}\left(\sum_{k=0}^{\infty}D_{k,\lambda}^{(r)}(x)S_{2,\lambda}(n,k)\right)\frac{t^{n}}{n!}.
\end{split}
\end{equation}
Therefore, by \eqref{36} and \eqref{37}, we obtain the following theorem.
\begin{thm}\label{Theorem 9}
For $n \ge 0$, we have
\begin{equation*}
\begin{split}
\beta_{n,\lambda}^{(r)}(x) = \sum_{k=0}^{n}D_{k,\lambda}^{(r)}(x)S_{2,\lambda}(n,k).
\end{split}
\end{equation*}
\end{thm}
By replacing $t$ by $\log_{\lambda}(1+t)$ in \eqref{08}, we get
\begin{equation}\label{38}
\begin{split}
\left(\frac{\log_{\lambda}(1+t)}{t}\right)^{r}(1+t)^{x} & = \sum_{k=0}^{\infty}\beta_{k,\lambda}^{(r)}(x)\frac{1}{k!}(\log_{\lambda}(1+t))^{k}\\
& = \sum_{k=0}^{\infty}\beta_{k,\lambda}^{(r)}(x)\sum_{n=k}^{\infty}S_{1,\lambda}(n,k)\frac{t^{n}}{n!}\\
& = \sum_{n=0}^{\infty}\left(\sum_{k=0}^{n}\beta_{k,\lambda}^{(r)}(x)S_{1,\lambda}(n,k)\right)\frac{t^{n}}{n!}.
\end{split}
\end{equation}
On the other hand,
\begin{equation}\label{39}
\begin{split}
\left(\frac{\log_{\lambda}(1+t)}{t}\right)^{r}(1+t)^{x} = \sum_{n=0}^{\infty}D_{n,\lambda}^{(r)}(x)\frac{t^{n}}{n!}.
\end{split}
\end{equation}
Therefore, by \eqref{38} and \eqref{39}, we obtain the following theorem
\begin{thm}\label{Theorem 10}
For $n \ge 0$, we have
\begin{equation*}
\begin{split}
D_{n,\lambda}^{(r)}(x) = \sum_{k=0}^{n}\beta_{k,\lambda}^{(r)}(x)S_{1,\lambda}(n,k).
\end{split}
\end{equation*}
\end{thm}
From \eqref{34}, we note that
\begin{equation}\label{40}
\begin{split}
\sum_{n=0}^{\infty}D_{n,\lambda}^{(r)}\frac{t^{n}}{n!}& = \underbrace{\left(\frac{\log_{\lambda}(1+t)}{t}\right)\times \cdots \times \frac{\log_{\lambda}(1+t)}{t}}_{r-\mbox{times}}\\
& = \sum_{n=0}^{\infty}\left(\sum_{l_{1}+\cdots +l_{r}= n}{n \choose l_{1},\cdots l_{r}}D_{l_{1},\lambda}\cdots D_{l_{r},\lambda}\right)\frac{t^{n}}{n!}.
\end{split}
\end{equation}

By \eqref{40}, we get
\begin{equation}\label{41}
\begin{split}
D_{n,\lambda}^{(r)} = \sum_{l_{1}+\cdots +l_{r}= n} {n \choose l_{1},\cdots l_{r}}D_{l_{1},\lambda}\cdots D_{l_{r},\lambda}, (n \ge 0).
\end{split}
\end{equation}

On the other hand, by \eqref{35}, we get
\begin{equation}\label{42}
\begin{split}
& \sum_{n=0}^{\infty}D_{n,\lambda}^{(r)}\frac{t^{n}}{n!} = \left(\frac{\log(1+t)}{t}\right)^{r}\int_{0}^{1}\cdots \int_{0}^{1}(1+t)^{\lambda(x_{1}+\cdots +x_{r})}dx_{1}\cdots dx_{r}\\
& \quad  = \left(\frac{\log(1+t)}{t}\right)^{r}\sum_{m=0}^{\infty}\lambda^{m}\frac{(\log(1+t))^{m}}{m!}\int_{0}^{1}\cdots\int_{0}^{1}(x_{1}+\cdots+x_{r})^{m}dx_{1}\cdots dx_{r}\\
& \quad=  \frac{1}{t^{r}}\sum_{m=0}^{\infty}\lambda^{m}\sum_{l_{1}+\cdots +l_{r}= m}{m \choose l_{1},\cdots ,l_{r}}\frac{1}{(l_{1}+1)\cdots(l_{r}+1)}\frac{(\log(1+t))^{m+r}}{m!}\\
& \quad=  \frac{1}{t^{r}}\sum_{m=0}^{\infty}\lambda^{m}\sum_{l_{1}+\cdots +l_{r}= m}{m \choose l_{1},\cdots ,l_{r}}\frac{1}{(l_{1}+1)\cdots(l_{r}+1)}\frac{(m+r)!}{m!}\\
& \quad \quad \quad \times \sum_{n=m+r}^{\infty}S_{1}(n,m+r)\frac{t^{n}}{n!} \\
& \quad=  \sum_{m=0}^{\infty}\lambda^{m}\sum_{l_{1}+\cdots +l_{r}= m}{m \choose l_{1},\cdots ,l_{r}}\frac{1}{(l_{1}+1)\cdots(l_{r}+1)}\frac{(m+r)!}{m!}\\
& \quad \quad \quad \times \sum_{n=m}^{\infty}S_{1}(n+r,m+r)\frac{t^{n}}{(n+r)!}\\
& \quad=  \sum_{n=0}^{\infty}\left(\sum_{m=0}^{n}\lambda^{m}\sum_{l_{1}+\cdots + l_{r}= m}\frac{{n \choose l_{1},\cdots ,l_{r}}}{(l_{1}+1)\cdots (1_{r}+1)}
S_{1}(n+r,m+r)\frac{{m+r\choose r}}{{n+r \choose r}}\right)\frac{t^{n}}{n!}.
\end{split}
\end{equation}
Therefore, by comparing the coefficients on both sides of \eqref{42}, we obtain the following theorem.
\begin{thm}\label{Theorem 11}
For $n \ge 0$, we have
\begin{equation*}
\begin{split}
&D_{n,\lambda}^{(r)} = \sum_{l_{1}+\cdots +l_{r}= m}{n \choose l_{1},\cdots ,l_{r}}D_{l_{1},\lambda}\cdots D_{l_{r},\lambda}\\
& = \quad \sum_{m=0}^{n}\lambda^{m}\sum_{l_{1}+\cdots +l_{r}= m}{m \choose l_{1},\cdots , l_{r}}\frac{S_{1}(n+r,m+r)}{(l_{1}+1)\cdots (1_{r}+1)}\frac{{m+r\choose r}}{{n+r \choose r}}.
\end{split}
\end{equation*}
\end{thm}

\medskip
\section{\bf Conclusion}

In the spirit of \cite{1}, we studied the degenerate Daehee polynomials and numbers which were recently introduced by Jang et. al \cite{4}. We derived their explicit expressions and some identities involving them. Further, we introduced the multiple degenerate Daehee numbers and higher-order degenerate Daehee polynomials and numbers and deduced their explicit expressions and some identities related to them. \\
\indent The possible applications of our results are as follows. The first one is their applications to identities of symmetry. For example,  in \cite{12} many symmetric identities in three variables, related to degenerate Euler polynomials and alternating generalized falling factorial sums, were obtained. The second one is their applications to differential equations from which we can derive some useful identities. For example, in \cite{9} an infinite family of nonlinear differential equations, having the generating function of the degenerate Bernoulli numbers as a solution, were derived. As a result, an identity, involving the degenerate Bernoulli and higher-order degenerate Bernoulli numbers, were obtained. Similar things had been done for the degenerate Euler numbers.
The third one is their applications to probability theory. Indeed,  in \cite{13} it was shown that both the degenerate $\lambda$-Stirling polynomials of the second and the $r$-truncated degenerate $\lambda$-Stirling polynomials of the second kind appear in certain expressions of the probability distributions of appropriate random variables. \\
\indent Finally, it is one of our future projects to continue to study various degenerate versions of some special polynomials and their applications to mathematics, science and engineering.

\vspace{6pt}
\section{ author contributions}
 T.K. and D.S.K. conceived of the framework and structured the whole paper; D.S.K. and T.K. wrote the paper; J.K. and H.Y.K. checked the results of the paper;
D.S.K. and T.K. completed the revision of the article. All authors have read and agreed to the published version of the manuscript.

\section{ funding}
 This work was supported by the National Research Foundation of Korea (NRF) grant funded by the
Korea government (MEST) (No. 2017R1E1A1A03070882).

\section{ conflicts of interest}
The authors declare that they have no competing interests.

\end{document}